\documentclass{article}

\usepackage{longtable}
\usepackage{amssymb, amsthm, amsmath}
\usepackage{xy} \input xy \xyoption{all}
 \usepackage{times}

\newtheorem{theorem}{Theorem}
\theoremstyle{definition}
\newtheorem*{algorithm}{Algorithm}
\newtheorem*{definition}{Definition}
\newtheorem*{remark}{Remark}
\newtheorem*{analysis}{Analysis}

\newcommand{\Z}{\mathbb{Z}}

\newcommand{\Q}{\mathbb{Q}}

\begin{document}

\title{Decomposing replicable functions}

\author{John M\raise0.4ex\hbox{\scriptsize\sc{c}}Kay\footnote{The first author is supported by NSERC.}\\Department of Computer Science\\Concordia University\\Montreal H3G 1M8, QC, Canada\\mckay@cs.concordia.ca \and David Sevilla\footnote{The second author is partially supported by Spanish Ministry of Science grant MTM2004-07086.}
\\Johann Radon Institute for Computational and Applied Mathematics\\Altenbergerstrasse 69\\A-4040 Linz, Austria\\david.sevilla@ricam.oeaw.ac.at}

\date{}

\maketitle

\begin{abstract}
We describe an algorithm to decompose rational functions from which we
determine the poset of groups fixing these functions.
\end{abstract}

\section{Introduction to replicable functions}

We assume familiarity with the notation and contents of \cite{ACMS92},
\cite{ConNor79}, \cite{ForMcKNor94} and \cite{McK01}. Replicable functions are
definable in terms of a generalized Hecke operator or by constraints on their
coefficients, see \cite{Nor84}. Each such function $f$ is fixed by a group
$G_f$ that is commensurable with the modular group, $\mathrm{PSL}(2,\Z)$.
There is a natural poset formed by these stabilizer groups, considered up to
conjugation by $z\mapsto kz,\ k\in\Z^{>0}$.

Replicable functions have a $q$-series (Fourier series) expansion at $i\infty$
of the form
\[f(q) = \frac 1 q + \sum_{k\geq1} a_k\,q^k\ , \qquad q = e^{2\pi iz},\ \Im(z) > 0, \quad \forall\ k:\ \ a_k\in\Z .\]
Cummins proves in \cite{Cum93} that a finite series implies that $f(q)=1/q +
cq, c\in\{0,1,-1\}$ -- the modular fictions, exp, cos, sin which we shall
hereafter ignore.

Computations suggest that there are 616 other replicable functions of which
171 are monstrous moonshine functions, see \cite{ConMcKSeb04}. There is no
satisfactory proof of the completeness of this list although it is compatible
with several independent computational checks.

The following remarkable result of Norton is fundamental, see \cite{Nor84},
\cite{Cum93}.

\begin{theorem}
A replicable function is determined by its coefficients in the Norton basis,
$\{a_k\}\,,\ k\in B=\{1,2,3,4,5,7,8,9,11,17,19,23\}$.
\end{theorem}

\section{Algorithms for functional decomposition and computation of relations}

\subsection{Functional decomposition}

We sketch the theory of univariate rational decomposition and an algorithm for
the computation of the poset of replicable functions with respect to rational
relations.

\begin{definition}
In $T=\Q(t)\setminus\Q$ we define the binary operation of \textbf{composition}
as
\[g(t)\circ h(t)=g(h(t))=g(h)(t).\]
$(T,\circ)$ is a semigroup with $t$ as neutral element.

If $f=g\circ h$, we call this a \textbf{decomposition} of $f$ and say that $g$
is a \textbf{left component} of $f$ and $h$ is a \textbf{right component} of
$f$. A decomposition is \textbf{trivial} if $g$ or $h$ is a unit with respect
to composition.

Decompositions $f=g_1\circ h_1=(g_1\circ u)\circ(u^{-1}\circ h_2)$ are called
\textbf{equivalent}, where $u$ is invertible with respect to composition and
$u^{-1}$ is its functional inverse.

Given a rational function $f\in T$, we call it \textbf{indecomposable} if it
is not a unit and all its decompositions are trivial. A decomposition of
$f\in\Q(t)$ of length $r$, $f=g_1\circ\cdots\circ g_r$, is called
\textbf{refined} if each $g_i$ is indecomposable.
\end{definition}

The units with respect to composition are linear fractional transformations.
The decomposition problem is: given $f\in\Q(t)$, compute all the
decompositions of $f$, i.e., find a representative $(h_i,g_i)$ for each class
of decompositions with respect to the equivalence relation above. Solving this
problem leads to the computation of all refined decompositions.

\begin{definition}
For a non-constant rational function $f(t)=f_N(t)/f_D(t)$ with
$f_N,f_D\in\Q[t]$ and $\gcd(f_N,f_D)=1$ we define the \textbf{degree} of $f$
as
\[\deg f=\max\{\deg f_N,\ \deg f_D\}.\]

We also define $\deg a=0$ for all non-zero $a\in\Q$.
\end{definition}

Because the solution to the problem may not be unique, most decomposition
algorithms have two steps: first, we compute candidates for the right
components, then check for their associated left components.

\begin{remark}
Given $f,h\in\Q(t)$, we can efficiently test if there is a $g\in\Q(t)$ with
$f=g\circ h$. It is necessary that $\deg h$ divides $\deg f$. We then solve
the equations resulting from the $q-$expansion of $f-(g\circ h)$. This is fast
as the equations are linear.
\end{remark}

We introduce a useful notion that will be the starting point for the
decomposition algorithm, see \cite{GutRubSev01} and \cite{Sev04}.

\begin{definition}
Let $f=f_N/f_D\in\Q(t)$ with $f_N,f_D\in\Q[t]\ ,\ \gcd(f_N,f_D)=1$. We say
that $f$ is in \textbf{normal form} when $\deg f_N>\deg f_D$ and $f_N(0)=0$
(or simply, $f(\infty)=\infty$, $f(0)=0$).
\end{definition}

\begin{theorem}
Let $f\in T$.
\begin{itemize}
\item[(i)]  There exist units $u,v\in\Q(t)$ such that $u\circ f\circ v$ is
in normal form, with both numerator and denominator monic.
\item[(ii)] Let $f\in\Q(t)$ be in normal form. If $f=g\circ h$, there is a
unit $u$ such that $g\circ u$ and
$u^{-1}\circ h$ are in normal form.
\end{itemize}
\end{theorem}

The following is the key to the decomposition algorithm.

\begin{theorem}
Let $f,g,h\in\Q(t)$ with $f=f_N/f_D$, $h=h_N/h_D$ where
$f_N,f_D,h_N,h_D\in\Q[t]$, $\gcd(f_N,f_D)=1$ and $\gcd(h_N,h_D)=1$, and
$f=g\circ h$. If $f,g,h$ are in normal form, then $h_N | f_N$ and $h_D | f_D$.
\end{theorem}

\begin{proof}
Let
\[g=\frac{t^r+c_{r-1}t^{r-1}+\cdots+c_1t}{d_{r-1}t^{r-1}+\cdots+d_0},\qquad
d_0\neq 0,\]
then
\[f=\frac{h_N^r+c_{r-1}h_N^{r-1}h_D+\cdots+c_1h_Nh_D^{r-1}}{d_{r-1}h_N^{r-1}h_D+\cdots+d_0h_D^r}\]
and, as the degree is multiplicative with respect to composition, there is no
simplification in this expression. The result follows.
\end{proof}

We describe the algorithm now.

\begin{algorithm}[Rational decomposition] $ $
\begin{description}
   \item [Input:] $f\in T$.
   \item [Output:] all non-trivial decompositions $(g,h)$ of $f$, if any exists.
\end{description}
\begin{description}
   \item[A] Compute $u$ and $v$ so that $\overline{f}=u\circ f \circ v$ is in
normal form. Let $f_N,f_D$ be the monic numerator and denominator of
$\overline{f}$.
   \item[B] Factor $f_N$ and $f_D$. From this compute
$D=\{(A_1,B_1),\ldots,(A_m,B_m)\}$, the set of pairs $(A,B)$ such that $A,B$
are monic polynomials dividing $f_N,f_D$ respectively. Set $i=1$.
   \item[C] Check if there exists $g\in\Q(t)$ with $\overline{f}=g(A_i/B_i)$;
if it does, add $\left(u^{-1}(g),h(v^{-1})\right)$ to the list of
decompositions of $f$.
   \item[D] If $i<m$, increase $i$ and go to \textbf{C}, otherwise
return the list of decompositions.
\end{description}
\end{algorithm}

\begin{analysis}
The description above shows that the algorithm correctly computes at least one
representative for each equivalence class of decompositions (an extra step
would be needed to avoid having more than one representative for each
decomposition class). The algorithm has exponential complexity due to the
possibility of having an exponential number of candidates in the worst case.
In practice, degree conditions reduce the number of candidates. In tests we
have found that about 85\% of the time is spent on the factoring, and the
number of candidates is small (random polynomials are irreducible). Because of
this, the algorithm is fast.
\end{analysis}

\subsection{Rational relations}

To find relations between two $q$-series we follow a simple procedure. We use
the fact that replicable functions correspond to groups acting on the upper
half plane, and a rational function of degree $n$ is an $n:1$ map.

\begin{algorithm}[Computation of rational relations] $ $
\begin{description}
\item [Input:] two $q$-series $s_1,s_2$ as described in the introduction.
\item [Output:] all rational relations of the form
$s_1(q^k)=f(s_2(q))$, $k\geq1$.
\end{description}
\begin{description}
\item [A] Compute the orders $e_1,\ldots,e_r$ of the generators
$M_1,\ldots,M_r$ of the fundamental region of $s_1$. The hyperbolic area of
the region is $A_1:=(r-2)\pi-\sum \pi/e_i$. Compute the area $A_2$ for $s_2$.
\item [B] If $d:=A_2/A_1$ is not an integer, then there are no
relations. Otherwise, put $r=1$.
\item [C] Let
\[f=\frac{t^d+a_{d-1}t^{d-1}+\cdots+a_0}{t^{d-r}+b_{d-r-1}t^{d-r-1}+\cdots+b_0}\]
and solve for $a_i,b_j$ the linear system given by $f(s_2(q^r))-s_1(q)$.
%
\item [D] If there is a non-trivial solution to the system, store the
corresponding $f$ and $r$. If $r<d$, increase $r$ and go to \textbf{C},
otherwise return all relations found.
\end{description}
\end{algorithm}

\begin{analysis}
In each relation, the degree of the numerator is the ratio of the areas, and
the difference between the degrees of numerator and denominator is the
exponent of $q$ in the function $s_2$. In step \textbf{A}, the orders of the
non-identity elements are determined by the trace squared divided by the
determinant. In step \textbf{C}, solving for the $a_i,b_j$ requires less than
$2d$ coefficients.
\end{analysis}

\begin{remark}
The values $k\geq1$ correspond to the conjugation $z\mapsto kz,\ k\in\Z^{>0}$,
that is, $q\mapsto q^k$.
\end{remark}

\section{The computation}

For each of the 616 replicable functions, we have:

\begin{itemize}
\item the coefficients $a_1,\ldots,a_{23}$,
\item parabolic and elliptic generators for the fixing groups, i.e. generators
of the stabilizers of the vertices of a fundamental region.
\end{itemize}

For each pair of series we determine whether there is a rational relation
between them as in Section 3.2. We decompose any rational relations as
described in Section 3.1 in order to refine the decompositions, we repeat
until we have all refined decompositions. In terms of the poset graph we use:

\begin{algorithm} [Poset refinement] $ $
\begin{description}
\item[A] Draw a vertex for each of the 616 functions.
\item[B] For each pair of functions $s_1$ and $s_2$, compute all
rational relations of the form $s_1(q^k)=f(s_2(q))$, if any. For each
relation, draw a labelled directed
edge
\[s_1\ \stackrel{\deg f,\;k}{\longrightarrow}\ s_2
\quad\mbox{where}\quad s_1(q^k)=f(s_2(q)).\]
\item[C] For each of these, compute all the decompositions
of $f$. For each decomposition
$(g,h)$ compute $s_3(q^j):=h(s_2(q))$, $j\geq1$ and replace the edge
with the two edges
\[s_1\ \stackrel{\deg g,\;k\!/\!j}{\longrightarrow}\ s_3\
\stackrel{\deg h,\;j}{\longrightarrow}\ s_2\]
\item[D] Repeat step \textbf{B} until all the rational functions are
indecomposable.
\end{description}
\end{algorithm}

The computations described in this section were performed in a Pentium-IV 2GHz
using Maple 7. First, once the areas were known, it was seen that the maximum
possible degree of a rational relation would be 96. The search for rational
relations took about 20 hours, about 10\% of this time was spent in
precomputing 200 coefficients for each series from the initial 23. We found
2419 relations in this step, we show their degrees below.
\[\begin{array}{cccc}
degree & number& degree & number\\
2 & 698 & 16 & 52 \\
3 & 243 & 18 & 60 \\
4 & 422 & 20 & 2 \\
5 & 26 & 24 & 71 \\
6 & 333 & 28 & 2 \\
8 & 178 & 30 & 8 \\
9 & 40 & 32 & 4 \\
10 & 14 & 36 & 40 \\
12 & 209 & 48 & 5 \\
14 & 4 & 72 & 2 \\
15 & 6 \\
\end{array}\]

In Step \textbf{C}, we decompose all the functions we found previously, remove
the repeated relations, and continue until all functions are indecomposable.
In this way, and since our decomposition algorithm outputs all possible
decompositions up to units, we ensure that we find all missing functions, if
any, from the lists available. The computation of all possible decompositions
is fundamental since there exists no formal proof of the completeness of our
initial data. Our computation provides this. The decomposition of all rational
relations took around 30 hours overall. In the end, we obtained 1049
indecomposable rational relations. We summarize them in Table \ref{table1} in
the appendix. We also list the connected components of the graph there.

For each function we give its immediate predecessors and successors. 
For each edge we give two numbers, the degrees of the numerator and
denominator of the rational relation; the first is the degree of the relation,
and the power of $q$ is given by the difference of the two numbers. For
example, $(1A,2:0)$ in line $2a$ means that $j(q^2)$ is a degree two
polynomial in the principal modulus $2a$. Notice that in some cases there are
two edges between two given functions.



\section{Remarks}

It is noteworthy that for two series $s_1,s_2$ we may find more than one
relation, i.e. $s_1(q^{k_1})=f_1(s_2(q))$ and $s_1(q^{k_2})=f_2(s_2(q))$ with
$k_1\neq k_2$. By computation of a resultant, we can find a polynomial
relation of the type $P(s_1(q),s_1(q^{k_1/k_2}))=0$.

Computation reveals a remarkable fact about the relation between $j$ (labelled
$f =1A$) and the principal modulus for $\Gamma(3)$, $t=s(z/3)$ where $s =
(\eta(q)/\eta(q^9))^3+3$, labelled $9B$. Specifically, refined decomposition
chains of different lengths exist for
%
%
%
%
\begin{figure}[ht!]
\[\xy
 (-80,18)*{f=\displaystyle\frac{t^3(t^3+6^3)^3}{(t^3-3^3)^3}\,;};
 (-65,9)*{\mbox{namely} \quad\displaystyle f=t^3\circ\frac{t(t-12)}{t-3}\circ\frac{t(t+6)}{t-3}};
 (-64,0)*{\mbox{and} \quad\displaystyle f=\frac{t^3(t+24)}{t-3}\circ\frac{t(t^2-6\,t+36)}{t^2+3\,t+9}\,.};
 {\ar@{-}^3 (0,0)*+{9B}; (20,6)*+{3B};};
 {\ar@{-}^4 (20,6)*+++{}; (0,18)*+{1A};};
 {\ar@{-}^3 (0,18)*+++{}; (-20,12)*+{3C};};
 {\ar@{-}^2 (-20,12)*+++{}; (-10,6)*+{9A};};
 {\ar@{-}^2 (-10,6)*+++{}; (0,0)*+++{};};
\endxy\]
\end{figure}

We believe this is the first example of a rational function in $\Q(t)$ with
refined decomposition chains of different lengths. This does not occur with
polynomials, see \cite{GutSev06}.

Norton points out that to every component (other than the fictions) there is
at least one function that is either monstrous or the translate of a monstrous
function.

\section*{Appendix}

\begin{longtable}{p{33ex}cp{33ex}}
\caption[Graph table extended]{Table of minimal relations sorted by middle function}  \\
\multicolumn{1}{c}{Up} & f & \multicolumn{1}{c}{Down}  \label{table1}\\
\\
\endfirsthead
\multicolumn{1}{c}{Up} & f & \multicolumn{1}{c}{Down}  \\
\\
\endhead
 &  $1A$  & (3:2,~$2B$), (3:1,~$2B$), (2:0,~$2a$), (4:3,~$3B$), (4:1,~$3B$), (3:0,~$3C$), (3:1,~$4A$), (6:5,~$5B$), (6:1,~$5B$), (5:0,~$5a$), (8:7,~$7B$), (8:1,~$7B$), (14:13,~$13B$), (14:1,~$13B$), (28:21,~$49a$),   \\
 &  $2A$  & (2:1,~$2B$), (2:0,~$4B$), (2:0,~$4a$), (4:3,~$6D$), (4:1,~$6D$), (3:0,~$6d$), (6:5,~$10C$), (6:1,~$10C$), (5:0,~$10b$),   \\
($1A$,~3:2), ($1A$,~3:1), ($2A$,~2:1),  &  $2B$  & (2:1,~$4C$), (2:0,~$4C$), (2:0,~$4D$), (4:3,~$6E$), (4:1,~$6E$), (3:0,~$6F$), (6:5,~$10E$), (6:1,~$10E$),   \\
($1A$,~2:0),  &  $2a$  & (3:2,~$4C$), (4:3,~$6c$), (4:1,~$6c$), (5:0,~$10c$), (8:7,~$14b$), (8:1,~$14b$),   \\
 &  $3A$  & (2:1,~$3B$), (3:2,~$6C$), (3:1,~$6C$), (2:0,~$6a$), (2:0,~$6b$), (3:0,~$9a$), (3:0,~$9b$), (3:1,~$12A$), (5:0,~$15b$), (8:7,~$21B$), (8:1,~$21B$),   \\
($1A$,~4:3), ($1A$,~4:1), ($3A$,~2:1),  &  $3B$  & (3:2,~$6E$), (3:1,~$6E$), (2:0,~$6c$), (3:2,~$9B$), (3:0,~$9B$), (3:0,~$9c$), (3:1,~$12B$),   \\
($1A$,~3:0),  &  $3C$  & (3:2,~$6F$), (3:1,~$6F$), (2:0,~$6\widetilde{a}$), (2:1,~$9A$), (3:0,~$9d$), (3:1,~$12D$), (6:5,~$15D$), (6:1,~$15D$),   \\
($1A$,~3:1), ($4\widetilde{a}$,~3:2), ($4\widetilde{b}$,~2:1),  &  $4A$  & (2:1,~$4C$), (2:0,~$8B$), (2:0,~$8a$), (4:3,~$12B$), (4:1,~$12B$), (3:0,~$12D$), (6:5,~$20C$), (6:1,~$20C$),   \\
\end{longtable}

The connected components of the graph are:

\noindent One with 480: $\{1A$, $2A$, $2B$, $2a$, $3A$, $3B$, $3C$, $4A$,
$4B$, $4C$, $4D$, $4a$, $4\widetilde{a}$, $4\widetilde{b}$, $5A$, $5B$, $5a$,
$6A$, $6B$, $6C$, $6D$, $6E$, $6F$, $6a$, $6b$, $6c$, $6d$, $6\widetilde{a}$,
$7A$, $7B$, $8A$, $8B$, $8C$, $8D$, $8E$, $8F$, $8a$, $8b$, $8c$,
$8\widetilde{a}$, $8\widetilde{b}$, $8\widetilde{c}$, $8\widetilde{d}$, $9A$,
$9B$, $9a$, $9b$, $9c$, $9d$, $10A$, $10B$, $10C$, $10D$, $10E$, $10a$, $10b$,
$10c$, $12A$, $12B$, $12C$, $12D$, $12E$, $12F$, $12G$, $12H$, $12I$, $12J$,
$12a$, $12b$, $12c$, $12d$, $12e$, $12f$, $12\widetilde{a}$,
$12\widetilde{b}$, $12\widetilde{c}$, $12\widetilde{d}$, $12\widetilde{e}$,
$12\widetilde{f}$, $12\widetilde{g}$, $12\widetilde{h}$, $12\widetilde{i}$,
$12\widetilde{j}$, $13A$, $13B$, $14A$, $14B$, $14C$, $14a$, $14b$, $14c$,
$15A$, $15B$, $15C$, $15D$, $15a$, $15b$, $16A$, $16B$, $16C$, $16a$, $16b$,
$16c$, $16d$, $16e$, $16f$, $16g$, $16h$, $16\widetilde{a}$,
$16\widetilde{b}$, $16\widetilde{c}$, $16\widetilde{d}$, $16\widetilde{e}$,
$18A$, $18B$, $18C$, $18D$, $18E$, $18a$, $18b$, $18c$, $18d$, $18e$, $18f$,
$18g$, $18h$, $18i$, $18j$, $18\widetilde{a}$, $20A$, $20B$, $20C$, $20D$,
$20E$, $20F$, $20a$, $20b$, $20c$, $20d$, $20e$, $20\widetilde{a}$,
$20\widetilde{b}$, $20\widetilde{c}$, $20\widetilde{d}$, $20\widetilde{e}$,
$20\widetilde{f}$, $20\widetilde{g}$, $20\widetilde{h}$, $21A$, $21B$, $21C$,
$21D$, $24A$, $24B$, $24C$, $24D$, $24E$, $24F$, $24G$, $24H$, $24I$, $24J$,
$24a$, $24b$, $24c$, $24d$, $24e$, $24f$, $24g$, $24h$, $24i$, $24j$,
$24\widetilde{a}$, $24\widetilde{b}$, $24\widetilde{c}$, $24\widetilde{d}$,
$24\widetilde{e}$, $24\widetilde{f}$, $24\widetilde{g}$, $24\widetilde{h}$,
$24\widetilde{i}$, $24\widetilde{j}$, $24\widetilde{k}$, $24\widetilde{l}$,
$24\widetilde{m}$, $24\widetilde{n}$, $24\widetilde{o}$, $24\widetilde{p}$,
$24\widetilde{q}$, $24\widetilde{r}$, $24\widetilde{s}$, $24\widetilde{t}$,
$25A$, $25a$, $26a$, $27A$, $27a$, $27b$, $27c$, $27d$, $27e$, $28A$, $28B$,
$28C$, $28D$, $28a$, $28\widetilde{a}$, $28\widetilde{b}$, $28\widetilde{c}$,
$28\widetilde{d}$, $30A$, $30B$, $30C$, $30D$, $30E$, $30F$, $30G$, $30a$,
$30b$, $30c$, $30d$, $30e$, $30f$, $30\widetilde{a}$, $32A$, $32B$, $32a$,
$32b$, $32c$, $32d$, $32e$, $32\widetilde{a}$, $32\widetilde{b}$,
$32\widetilde{c}$, $32\widetilde{d}$, $32\widetilde{e}$, $35a$, $36A$, $36B$,
$36C$, $36D$, $36a$, $36b$, $36c$, $36d$, $36e$, $36f$, $36g$, $36h$, $36i$,
$36\widetilde{a}$, $36\widetilde{b}$, $36\widetilde{c}$, $36\widetilde{d}$,
$36\widetilde{e}$, $36\widetilde{f}$, $36\widetilde{g}$, $36\widetilde{h}$,
$36\widetilde{i}$, $36\widetilde{j}$, $36\widetilde{k}$, $36\widetilde{l}$,
$36\widetilde{m}$, $36\widetilde{n}$, $36\widetilde{o}$, $36\widetilde{p}$,
$36\widetilde{q}$, $36\widetilde{r}$, $36\widetilde{s}$, $39B$, $40A$, $40B$,
$40C$, $40a$, $40b$, $40c$, $40d$, $40e$, $40\widetilde{a}$,
$40\widetilde{b}$, $40\widetilde{c}$, $40\widetilde{d}$, $40\widetilde{e}$,
$40\widetilde{f}$, $40\widetilde{g}$, $40\widetilde{h}$, $40\widetilde{i}$,
$40\widetilde{j}$, $40\widetilde{k}$, $42C$, $42a$, $42b$, $42c$, $42d$,
$45A$, $45a$, $45b$, $45c$, $48A$, $48a$, $48b$, $48c$, $48d$, $48e$, $48f$,
$48g$, $48h$, $48\widetilde{a}$, $48\widetilde{b}$, $48\widetilde{c}$,
$48\widetilde{d}$, $48\widetilde{e}$, $48\widetilde{f}$, $48\widetilde{g}$,
$48\widetilde{h}$, $48\widetilde{i}$, $48\widetilde{j}$, $48\widetilde{k}$,
$48\widetilde{l}$, $48\widetilde{m}$, $48\widetilde{n}$, $48\widetilde{o}$,
$49a$, $50A$, $50a$, $52\widetilde{a}$, $52\widetilde{b}$, $54A$, $54a$,
$54b$, $54c$, $54d$, $54\widetilde{a}$, $56A$, $56B$, $56a$, $56b$, $56c$,
$56\widetilde{a}$, $56\widetilde{b}$, $56\widetilde{c}$, $56\widetilde{d}$,
$56\widetilde{e}$, $56\widetilde{f}$, $56\widetilde{g}$, $60A$, $60B$, $60C$,
$60D$, $60E$, \linebreak $60F$, $60a$, $60b$, $60c$, $60d$, $60e$,
$60\widetilde{a}$, $60\widetilde{b}$, $60\widetilde{c}$, $60\widetilde{d}$,
$60\widetilde{e}$, $60\widetilde{f}$, $60\widetilde{g}$, $60\widetilde{h}$,
$60\widetilde{i}$, $60\widetilde{j}$, $60\widetilde{k}$, \linebreak
$60\widetilde{l}$, $60\widetilde{m}$, $60\widetilde{n}$, $63a$,
$63\widetilde{a}$, $64a$, $72a$, $72b$, $72c$, $72d$, $72e$,
$72\widetilde{a}$, $72\widetilde{b}$, $72\widetilde{c}$, $72\widetilde{d}$,
$72\widetilde{e}$, $72\widetilde{f}$, $72\widetilde{g}$, $72\widetilde{h}$,
$72\widetilde{i}$, $72\widetilde{j}$, $72\widetilde{k}$, $72\widetilde{l}$,
$72\widetilde{m}$, $72\widetilde{n}$, $72\widetilde{o}$, $72\widetilde{p}$,
$72\widetilde{q}$, $72\widetilde{r}$, $72\widetilde{s}$, $72\widetilde{t}$,
$80a$, $80\widetilde{a}$, $80\widetilde{b}$, \linebreak $80\widetilde{c}$,
$80\widetilde{d}$, $80\widetilde{e}$, $80\widetilde{f}$, $84C$,
$84\widetilde{a}$, $84\widetilde{b}$, $84\widetilde{c}$, $84\widetilde{d}$,
$84\widetilde{h}$, $90a$, $90b$, $90\widetilde{a}$, $96a$, $96\widetilde{a}$,
\linebreak $96\widetilde{b}$, $100\widetilde{a}$, $100\widetilde{b}$,
$100\widetilde{c}$, $100\widetilde{d}$, $104\widetilde{a}$,
$108\widetilde{a}$, $108\widetilde{b}$, $108\widetilde{c}$,
$108\widetilde{d}$, $108\widetilde{e}$, $108\widetilde{f}$,
$108\widetilde{g}$, $108\widetilde{h}$, $108\widetilde{i}$,
$108\widetilde{j}$, $108\widetilde{k}$, $112\widetilde{a}$,
$112\widetilde{b}$, $112\widetilde{c}$, $120a$, $120\widetilde{a}$,
$120\widetilde{b}$, $120\widetilde{c}$, $120\widetilde{d}$,
$120\widetilde{e}$, $120\widetilde{f}$, $120\widetilde{g}$,
$120\widetilde{h}$, $120\widetilde{i}$, $120\widetilde{j}$,
$120\widetilde{k}$, $120\widetilde{l}$, $120\widetilde{m}$,
$120\widetilde{n}$, $126a$, $126\widetilde{a}$, $140\widetilde{e}$,
$144\widetilde{a}$, $144\widetilde{b}$, $144\widetilde{c}$, \linebreak
$144\widetilde{d}$, $144\widetilde{e}$, $156\widetilde{b}$,
$160\widetilde{a}$, $160\widetilde{b}$, $168\widetilde{a}$,
$168\widetilde{b}$, $168\widetilde{e}$, $180\widetilde{a}$,
$180\widetilde{b}$, $180\widetilde{c}$, $180\widetilde{d}$,
$180\widetilde{e}$, $180\widetilde{f}$, $180\widetilde{g}$,
$196\widetilde{a}$, $216\widetilde{a}$, $216\widetilde{b}$,
$216\widetilde{c}$, $252\widetilde{a}$, $252\widetilde{b}$,
$252\widetilde{c}$, $504\widetilde{a}\}$.\vspace{1ex}

\noindent One with 20: $\{11A$, $22A$, $22B$, $22a$, $33A$, $33B$, $44A$,
$44a$, $44b$, $44c$, $44\widetilde{a}$, $44\widetilde{b}$, $66a$, $88A$,
$88\widetilde{a}$, $88\widetilde{b}$, $88\widetilde{c}$, $132\widetilde{a}$,
$132\widetilde{b}$, $264\widetilde{a}\}$.\vspace{1ex}

\noindent Five with 6: $\{23A$, $46A$, $46C$, $92A$, $92\widetilde{a}$,
$92\widetilde{b}\}$, $\{26A$, $26B$, $52A$, $52B$, $52a$, $104A\}$, $\{42A$,
$42B$, $42D$, $84A$, $84B$, $84a\}$, $\{52\widetilde{c}$, $52\widetilde{d}$,
$104\widetilde{b}$, $104\widetilde{c}$, $104\widetilde{d}$,
$208\widetilde{a}\}$, $\{84\widetilde{e}$, $84\widetilde{f}$,
$84\widetilde{g}$, $168\widetilde{c}$, $168\widetilde{d}$,
$168\widetilde{f}\}$. \vspace{1ex}

\noindent Ten with 3: $\{19A$, $38a$, $57A\}$, $\{35A$, $35B$, $70a\}$,
$\{39A$, $39C$, $117a\}$, $\{66A$, $66B$, $132a\}$, $\{70A$, $70B$, $140a\}$,
$\{76\widetilde{a}$, $152\widetilde{a}$, $228\widetilde{a}\}$,
$\{132\widetilde{c}$, $132\widetilde{d}$, $264\widetilde{b}\}$,
$\{140\widetilde{a}$, $140\widetilde{c}$, $280\widetilde{a}\}$,
$\{140\widetilde{b}$, $140\widetilde{d}$, $280\widetilde{b}\}$,
$\{156\widetilde{a}$, $156\widetilde{c}$, $468\widetilde{a}\}$.\vspace{1ex}

\noindent Sixteen with 2: $\{17A$, $34a\}$, $\{29A$, $58a\}$, $\{31A$,
$93A\}$, $\{34A$, $68A\}$, $\{38A$, $76a\}$, $\{41A$, $82a\}$, $\{51A$,
$102a\}$, $\{68\widetilde{a}$, $136\widetilde{a}\}$, $\{68\widetilde{b}$,
$136\widetilde{b}\}$, $\{76\widetilde{b}$, $152\widetilde{b}\}$, $\{78A$,
$78B\}$, $\{116\widetilde{a}$, $232\widetilde{a}\}$, $\{124\widetilde{a}$,
$372\widetilde{a}\}$, $\{156\widetilde{d}$, $156\widetilde{e}\}$,
$\{164\widetilde{a}$, $328\widetilde{a}\}$, $\{204\widetilde{a}$,
$408\widetilde{a}\}$.\vspace{1ex}

\noindent 24 isolated: $47A$, $55A$, $59A$, $62A$, $69A$, $71A$, $87A$, $94A$,
$95A$, $105A$, $110A$, $119A$, $124\widetilde{b}$, $188\widetilde{a}$,
$188\widetilde{b}$, $220\widetilde{a}$, $220\widetilde{b}$,
$236\widetilde{a}$, $276\widetilde{a}$, $284\widetilde{a}$,
$348\widetilde{a}$, $380\widetilde{a}$, $420\widetilde{a}$,
$476\widetilde{a}$.


\section*{Acknowledgement}

We thank Simon P. Norton for his help and insight.

\bibliographystyle{jcm}

\begin{thebibliography}{99}

\bibitem{ACMS92} \textsc{D. Alexander}, \textsc{C. J. Cummins}, \textsc{J. McKay} and \textsc{C.
Simons}, `Completely replicable functions', \emph{Groups, combinatorics \&
geometry}, London Math. Soc. Lecture Note Ser. 165 (Cambridge Univ. Press,
Durham, 1992), 87--98.


\bibitem{ConMcKSeb04} \textsc{J. H. Conway}, \textsc{J. McKay} and \textsc{A. Sebbar}, `On
the discrete groups of moonshine', \emph{Proc. Amer. Math. Soc.} 132 no. 8
(2004), 2233--2240.

\bibitem{ConNor79} \textsc{J. H. Conway} and \textsc{S. P. Norton}, `Monstrous
moonshine', \emph{Bull. London Math. Soc.} 11 (1979), 308--339.


\bibitem{Cum93} \textsc{C. J. Cummins}, `Some comments on replicable
functions', \emph{Modern trends in Lie algebra representation theory} (Queen's
Univ., Kingston, ON, 1994) 48--55, \emph{Queen's Papers in Pure and Appl.
Math.} 94 (1994).

\bibitem{ForMcKNor94} \textsc{D. Ford}, \textsc{J. McKay} and \textsc{S. P. Norton}, `More on
replicable functions'. \emph{Comm. in Algebra} 22 (1994) 5175-5193.

\bibitem{GutRubSev01} \textsc{J. Gutierrez}, \textsc{R. Rubio} and \textsc{D. Sevilla},
`Unirational fields of transcendence degree one and functional decomposition',
\emph{Proceedings of International Symposium on Symbolic and Algebraic
Computation}, ISSAC 2001, 167--175.

\bibitem{McK01} \textsc{J. McKay}, `Essentials of monstrous moonshine',
\emph{Groups and combinatorics---in memory of Michio Suzuki, Adv. Std. Pure
Math.} 32 (2001) 347--353.

\bibitem{Nor84} \textsc{S. P. Norton}, `More on moonshine',
\emph{Computational group theory} (London Academic Press, 1984) 185--193.

\bibitem{GutSev06} \textsc{J. Gutierrez} and \textsc{D. Sevilla}, `Building
counterexamples to generalizations for rational functions of Ritt's
decomposition Theorem', \emph{Journal of Algebra} 303 no. 2 (2006) 655--667.

\bibitem{Sev04} \textsc{D. Sevilla}, `Teoremas de Ritt y computaci\'on de
cuerpos unirracionales' (Ph. D. Thesis, University of Cantabria, 2004).

\end{thebibliography}

\end{document}